\documentclass{amsart}
\usepackage{amsmath,amsfonts,amscd,amssymb,epsf,mathrsfs}

\newtheorem*{theorem*}{Theorem}
\newtheorem{lemma}{Lemma}[section]

\theoremstyle{definition}

\theoremstyle{remark}
\newtheorem{remark}[lemma]{Remark}
\numberwithin{equation}{section}

\providecommand{\abs}[1]{\left\lvert#1\right\rvert}
\providecommand{\tabs}[1]{\lvert#1\rvert}

\newcommand{\ud}{\mathrm{d}}

\newcommand{\field}[1]{\ensuremath{\mathbb{#1}}}
\newcommand{\CC}{\field{C}}

\newcommand{\RR}{\field{R}}

\DeclareMathOperator{\im}{Im}

 \DeclareMathOperator{\re}{Re}

\newcommand{\del}{\partial}

\newcommand{\Ga}{\Gamma}

\newcommand{\hp}{\mathbb{H}}
\newcommand{\R}{{\mathbb{R}}}
\newcommand{\C}{{\mathbb{C}}}
\newcommand{\bk}{\backslash}
\newcommand{\pa}{\partial}
\newcommand{\g}{\gamma}

\newcommand{\ov}{\overline}

\newcommand{\vep}{\varepsilon}
\newcommand{\z}{\ov{z}}
\newcommand{\zbar}{\ov{z}}
\newcommand{\wbar}{\ov{w}}

\newcommand{\qf}{QF}
\newcommand{\teich}{T}

\newcommand{\univqf}{\mathscr{QF}}

\newcommand{\autoform}{B}
\newcommand{\holoautoform}{A}
\newcommand{\harmform}{\Omega}

\newcommand{\canclass}{\omega}
\newcommand{\delbar}{\ov{\pa}}

\providecommand{\define}[1]{\emph{#1}}

\newcommand{\n}{$n$\nobreakdash-\hspace{0pt}} 

\begin{document}

\title[Holomorphic factorization of determinants of Laplacians]
{Holomorphic factorization of determinants of Laplacians using
quasi-Fuchsian uniformization}
\author{Andrew Mcintyre}
\address{Centre de Recherches Math$\acute{\text{e}}$matiques
Universit$\acute{\text{e}}$ de Montr$\acute{\text{e}}$al, Pavillon
Andr$\acute{\text{e}}$-Aisenstadt, 2920 Chemin de la tour,
Montr$\acute{\text{e}}$al, PQ, H3T 1J4 Canada
}\email{mcintyre@crm.umontreal.ca}
\author{Lee-Peng Teo}\address{Faculty of Information Technology,
Multimedia University, Jalan Multimedia, Cyberjaya, 63100,
Selangor Darul Ehsan, Malaysia}\email{lpteo@mmu.edu.my}
\date{\today}
\subjclass[2000]{  30F60, 30F10, 30F30} \keywords{holomorphic
factorization, Laplacian, period matrix, differentials,
quasi-Fuchsian}
\begin{abstract} For a quasi-Fuchsian group $\Ga$ with ordinary
set $\Omega$, and $\Delta_{n}$ the
Laplacian on \n differentials on $\Ga\bk\Omega$,
we define a notion of a Bers dual basis $\phi_{1},\dotsc,\phi_{2d}$ for
$\ker\Delta_{n}$. We prove that
$\det\Delta_{n}/\det \langle\phi_{j},\phi_{k}\rangle$, is, up to an anomaly computed
by Takhtajan and the second author
in \cite{TT1}, the modulus squared of a holomorphic function
$F(n)$, where $F(n)$ is a quasi-Fuchsian analogue of the Selberg zeta
$Z(n)$.   This generalizes the D'Hoker-Phong formula
$\det\Delta_{n}=c_{g,n}Z(n)$, and is a
quasi-Fuchsian counterpart of the result for Schottky groups
proved
by Takhtajan and the first author
in \cite{MT}.
\end{abstract}

\maketitle

\section{Introduction}

Let $\Ga$ be a Fuchsian group, such that $\Ga\bk\hp_+\simeq X$
is a compact Riemann surface of genus $g>1$, let $X$ be
endowed with the hyperbolic metric, and let
$\Delta_{n}$ be the Laplacian acting on the $n$th power of
the canonical bundle of $X$, $n>1$.  (See Section \ref{prelim}
for definitions.)  In \cite{DhokerPhong},
D'Hoker and Phong computed the zeta-regularized determinant
of $\det\Delta_{n}$ in terms of the group $\Ga$:
\begin{equation}\label{DPint}
\det\Delta_{n} = c_{g,n}Z(n)
= c_{g,n}\prod_{\{\g\}} \prod_{m=0}^{\infty} ( 1 -\lambda(\g)^{n+m}),
\end{equation}
where $\{\g\}$ runs over primitive conjugacy classes in $\Ga$,
$\lambda(\g)$ is the multiplier of $\g$, and $c_{g,n}$ is a known
constant.  This theorem is an analogue of the expression, (essentially due
to Kronecker), when $X$ has
genus $1$ and a flat metric of area $1$,
\begin{equation}\label{Kronint}
{\det}'\Delta_{0}(\tau)
=4\im\tau\abs{\eta(\tau)}^{4},
\end{equation}
where $X=\langle z\mapsto z+1, z\mapsto z+\tau\rangle\bk\C$,
and $\eta$ is the Dedekind function,
\[
\eta(\tau)=q^{\tfrac{1}{24}}
 \prod_{m=1}^{\infty}(1-q^{m}),\quad q=e^{2\pi i\tau}.
\]

The expression \eqref{DPint} is a higher genus analogue of \eqref{Kronint}, but
other analogues --- in some sense closer to \eqref{Kronint}
--- are possible.

In \cite{ZT}, Takhtajan and Zograf
proved the following formula on the Teichm\"uller space $\teich_{g}$:
\begin{equation}\label{lo}
\pa_{\mu}\pa_{\ov{\nu}}\log\frac{\det\Delta_n}{\det N_{n}}
=\frac{6n^2-6n+1}{12\pi}
\langle\mu,\nu\rangle_{\scriptscriptstyle{WP}}.
\end{equation}
Here $\mu$ and $\nu$ are holomorphic tangent vectors,
$(N_n)_{jk}=\langle\phi_j,\phi_k\rangle$ for
a holomorphically varying basis $\phi_1,\dotsc,\phi_d$ of
holomorphic \n differentials, ($d=(2n-1)(g-1)$), which we call the period matrix, and
$\langle \,\cdot\;,\;\cdot\,\rangle_{\scriptscriptstyle{WP}}$
is the Weil-Petersson metric on
$\teich_{g}$.

Now, for each point $[X]\in\teich_{g}$, there is a unique normalized
marked Schottky group $\Ga$ with ordinary set $\Omega$, such that
$X\simeq\Ga\bk\Omega$. As a corollary to their construction of the
Liouville action functional for Schottky groups, Takhtajan and
Zograf constructed the classical Liouville action \cite{ZT87b}, a
positive real-valued function $S$ on the Schottky space (whose
points correspond to marked Schottky groups), and hence on the
Teichm\"uller space $\teich_{g}$ as well. The negative of this
function is a potential for the Weil-Petersson metric, namely,
\begin{align*}
\pa_{\mu}\pa_{\ov{\nu}}(-S)
=\langle\mu,\nu\rangle_{\scriptscriptstyle{WP}}.
\end{align*}
Therefore,
\begin{align*}
\log\frac{\det\Delta_n}{\det N_n} +\frac{6n^2-6n+1}{12\pi}S
\end{align*}
is a pluriharmonic real-valued function on $\teich_{g}$. Since
$\teich_{g}$ is contractible,
there exists a holomorphic function
$F(n)$ on $\teich_{g}$ such that
\[
\log\frac{\det\Delta_n}{\det N_n} +\frac{6n^2-6n+1}{12\pi}S=\log\abs{F(n)}^2,
\]
or equivalently,
\[
\frac{\det\Delta_n}{\det N_n}
=\exp{\left(-\frac{6n^2-6n+1}{12\pi}S\right)}\abs{F(n)}^2,
\]
which is what we refer to as ``holomorphic factorization''.

In \cite{MT}, L.~Takhtajan
and the first author showed that, with a suitable explicit choice of the basis
$\phi_1,\dotsc,\phi_d$, depending on the Schottky group
$\Ga$, the function $F(n)$ has the expression
\[
\begin{split}
F(n)&=(1-\lambda(L_{1}))^2(1-\lambda(L_{1})^2)^2
          \dotsb(1-\lambda(L_{1})^{n-1})^2(1-\lambda(L_{2})^{n-1})\times\\
&\quad\quad\cdot\prod_{\{\gamma\}}\prod_{m=0}^{\infty}(1-\lambda(\gamma)^{n+m}),
\end{split}
\]
where $\{\g\}$ runs over primitive conjugacy classes in the Schottky
group $\Ga$, and $L_1,\dotsc,L_g$ are generators
of $\Ga$ corresponding to the marking.  (The factors before the product sign are related to
the explicit choice of basis $\phi_1,\dotsc,\phi_d$.)
This theorem is closer in form to the genus $1$ result
than is \eqref{DPint},
and in fact specializes to it when $g=1$.

In \cite{TT1}, L.~Takhtajan and the second author extended the
results of \cite{ZT87b} to quasi-Fuchsian groups.  In particular,
they constructed the classical Liouville action,
a positive real-valued function $S$
on the quasi-Fuchsian space $\qf_{g}$ (whose points correspond to
marked quasi-Fuchsian groups).  We have $\qf_{g}\simeq\teich_{g}\times\teich_{g}$.
As in the Schottky case, the negative of this function is a
potential for the Weil-Petersson metric.
Therefore one may expect that a result similar to \cite{MT} will hold
for quasi-Fuchsian groups.  In this paper we prove such a result.

First it is necessary to make an appropriate choice of basis for
holomorphic \n differentials.  For this purpose we define the notion of
a \emph{Bers dual} basis.
Let $\Ga$ be a cocompact quasi-Fuchsian group with domain of
discontinuity $\Omega=\Omega_+\sqcup\Omega_-$, simultaneously
uniformizing the compact Riemann surfaces $X_\pm\simeq \Ga\bk\Omega_\pm$.
For $n>1$,  Bers \cite{Bers2}
introduced an invertible integral operator $K_-$ which maps the conjugate of
a holomorphic \n differential on $X_-$ to a holomorphic
\n differential on $X_+$. Its kernel is given by
\[
K_-(z,w)
= \frac{2^{2n-2}(2n-1)}{\pi}\sum_{\gamma\in \Ga}
\frac{\gamma'(z)^n}{(\gamma z-w)^{2n}},\quad\quad z\in \Omega_+, w\in \Omega_-.
\]
Given a basis  $\phi^+_1,\ldots,\phi^+_d$ of holomorphic \n
differentials of $X_+\simeq\Ga\bk\Omega_+$, the Bers dual is the
basis $\phi^-_1,\ldots, \phi^-_d$ of holomorphic \n differentials of
$X_-\simeq\Ga\bk\Omega_-$ such that
\[
K_-(z,w)=\sum_{k=1}^d \phi^+_k(z)\phi^-_k(w).
\]

Now, the Local Index
Theorem \eqref{lo} and the result in \cite{TT1} imply that on the
quasi-Fuchsian space, there exists a holomorphic
function $F(n)$ such that
\begin{align*}
\frac{\det\Delta_n(X_+)}{\det N_{n}([X_+])}\frac{\det\Delta_n(X_-)}{\det N_{n}([X_-])}
=\exp\left(-\frac{6n^2-6n+1}{12\pi}S\right)\abs{F(n)}^2.
\end{align*}
Here $(N([X_\pm]))_{jk}=\langle\phi^\pm_j,\phi^\pm_k\rangle$ for
Bers dual bases $\phi^+_1,\ldots,\phi^+_d$ and
$\phi^-_1,\ldots,\phi^-_d$. In this paper, we show that up to a
known multiplicative constant, the function $F(n)$ is given by the
product
\[
F(n) = \prod_{\{\g\}} \prod_{m=0}^{\infty} (1 - \lambda(\g)^{n+m}),
\]
where $\{\g\}$ runs over the set of conjugacy classes of
primitive elements of the quasi-Fuchsian group $\Ga$.

This theorem specializes to the D'Hoker-Phong result when restricted
to the  subspace of Fuchsian groups.  In fact, the method of proof
is to use \eqref{DPint}, together with pluriharmonicity (coming from
the Local Index Theorem), and symmetry properties of the various
quantities in the theorem under complex conjugation of the group
$\Ga\mapsto\ov{\Ga}$.

We collect the background facts we will need in Section
\ref{prelim}.  In following sections, we discuss the generalization
of the Local Index Theorem, the function $F(n)$, the Bers integral
operator and Bers dual basis.  The main theorem is
stated and proved in Section \ref{mainthm}.

\smallskip\noindent
\textbf{Acknowledgements.} We would like to thank Leon Takhtajan for
helpful suggestions.  A. Mcintyre would like to thank Paul Gauthier
for a useful discussion. The work of L.-P. Teo was partially
supported by MMU internal funding PR/2006/0590.

\section{Preliminaries}\label{prelim}

In this section, we collect some necessary facts and definitions. We
refer the reader to the references cited for further details.


\subsection{Quasi-Fuchsian groups and simultaneous uniformization}
\label{QFgroup}
(See \cite{Ahlfors2, Bers3, Bers4,Bers5, Kra2, TT1}.)
By definition, a \define{Kleinian group}
is a discrete subgroup $\Ga$ of the group of M\"obius transformations
$\mathrm{PSL}(2,\C)$
which acts properly discontinuously on some non-empty
open subset of the Riemann sphere $\widehat{\C}=\C\cup\{\infty\}$.
The largest such subset $\Omega\subset\widehat{\C}$
is called the \define{ordinary set} of $\Ga$ and its complement $\Lambda$
is called the \define{limit set} of $\Ga$.

A Kleinian group $\Ga$
is called a \define{quasi-Fuchsian group} if it leaves some directed Jordan curve
 $\mathcal{C}\in\widehat{\C}$
invariant.  We have $\Lambda\subset\mathcal{C}$; if
$\Lambda=\mathcal{C}$, the group is said to be \define{of the first kind}.
In this case, $\Omega$ consists of the two domains
$\Omega_{+}$ and $\Omega_{-}$ complementary to $\mathcal{C}$,
chosen such that the boundary of $\Omega_{\pm}$ (with the orientation
from the complex plane) is $\pm\mathcal{C}$.
If $\mathcal{C}$ is a circle or line, the group is called \define{Fuchsian}; in this
case, we will assume $\mathcal{C}$ has been
conjugated to $\widehat{\R}=\R\cup\{\infty\}$ with the usual orientation,
so $\Omega_{\pm}=\hp_{\pm}$, the upper (resp. lower) half plane.

It is known that a group $\Ga$ is quasi-Fuchsian if and only if it
is the deformation of a Fuchsian group under a quasiconformal
(q.~c.) map, that is,
\[
\Ga=w\Ga_{0} w^{-1}
\]
for some Fuchsian group $\Ga_{0}$ and some q.~c. homeomorphism
$w:\widehat{\C}\to\widehat{\C}$. (See the references for the
definition of a q.~c. map.)

An element $\gamma$ of $\mathrm{PSL}(2,\C)$ is called
\define{loxodromic} if it is conjugate in $\mathrm{PSL}(2,\C)$ to
$z\mapsto\lambda(\gamma)z$ for some $\lambda(\gamma)\in\C$ (called
the \define{multiplier}) such that $0<\tabs{\lambda(\gamma)}<1$; if
every element of a Kleinian group is loxodromic except the identity,
the group is called \define{totally loxodromic}.  If $\Ga$ is a
finitely generated, totally loxodromic quasi-Fuchsian group of the
first kind, then the quotient $\Ga\bk\Omega\simeq X$ has two
connected components $\Ga\bk\Omega_{\pm}\simeq X_{\pm}$, which are
nonsingular compact Riemann surfaces of common genus $g>1$.  If
$\Ga$ is Fuchsian as well, then $X_{-}\simeq\ov{X_{+}}$, the mirror
image of $X_{+}$ (obtained by taking the complex conjugate of all
local coordinates). Conversely, Bers simultaneous uniformization
theorem states that, given any two compact Riemann surfaces $X_{+}$,
$X_{-}$ of common genus $g>1$, there exists a finitely generated
totally loxodromic quasi-Fuchsian group of the first kind $\Ga$ such
that $\Ga\bk\Omega_{\pm}\simeq X_{\pm}$.

\emph{In the sequel, ``quasi-Fuchsian group'' will always refer to a finitely generated,
totally loxodromic quasi-Fuchsian group of the first kind.}

A \define{marked} quasi-Fuchsian group is a quasi-Fuchsian group $\Ga$
(with the above convention),
together with a choice $a_{1},\dotsc,a_{g},b_{1},\dotsc,b_{g}\in\Ga$
of generators corresponding to standard generators of $\pi_{1}(X_{+})$.
The marking directs $\mathcal{C}$ by taking the segment from the attracting
fixed point of $a_{1}$ to that of $a_{2}$ to be positively directed.
A marked quasi-Fuchsian group is \define{normalized} if the attracting
fixed points of $a_{1}$ and $a_{2}$ are $0$ and $1$ respectively, and the
repelling fixed point of $a_{1}$ is $\infty$.  A marked quasi-Fuchsian group
may be normalized by a (unique) overall conjugation in $\mathrm{PSL}(2,\C)$.


\subsection{\n differentials and determinants of Laplacians}
\label{detLaplacian}
(See \cite{TZ91} and references therein.)
Suppose that $\Ga$ is a quasi-Fuchsian group, with
$\Ga\bk\Omega_{\pm}\simeq X_{\pm}$
compact Riemann surfaces of genus $g>1$.
For integers $n$ and $m$, an \define{automorphic form of type $(n,m)$} is
a function $\phi:\Omega_{\pm}\to\widehat{\C}$ such that
\[
\phi(z)=\phi(\gamma z)\,\gamma'(z)^{n}
\overline{\gamma'(z)}^{m}\quad\quad
\text{for all}\quad z\in\Omega_{\pm},\;\gamma\in\Ga.
\]
We write $\autoform_{n,m}(\Omega_{\pm},\Ga)$, for the
space of smooth automorphic forms of type $(n,m)$, and
identify $\autoform_{n,m}(\Omega_{\pm},\Ga)$
with $\autoform_{n,m}(X_{\pm})$, the space of smooth sections of
$\canclass_{X_{\pm}}^n\otimes\overline{\canclass}_{X_{\pm}}^m$,
where $\canclass_{X_{\pm}}=T^{\ast}X_{\pm}$ is the holomorphic cotangent
bundle of $X_{\pm}$.
We abbreviate
$\autoform_{n,0}(\Omega_{\pm},\Ga)=\autoform_n(\Omega_{\pm},\Ga)$,
called \define{\n differentials}, and write
$\holoautoform_{n}(\Omega_{\pm},\Ga)$ for the
holomorphic \n differentials.

The hyperbolic metric on $X_{\pm}$, written locally as $\rho(z)\abs{\ud z}^{2}$,
induces a Hermitian metric
\begin{equation} \label{metric-n-diff}
\langle\phi,\psi\rangle =\int\limits_\mathcal{F}
\phi\overline{\psi}\rho^{1-n-m}\,\ud^{2}z,
\end{equation}
on $\autoform_{n,m}(X_{\pm})$, where $\mathcal{F}_{\pm}$ is a
fundamental region for $\Ga$ in $\Omega_{\pm}$, and
$\ud^{2}z=\frac{i}{2}\ud z\wedge\ud\zbar$ is the Euclidean area form
on $\Omega_{\pm}$.  The complex structure and metric determine a
connection
\[
D=\del_n\oplus\delbar_n:
\autoform_n(X_{\pm})\to
\autoform_{n+1}(X_{\pm})\oplus\autoform_{n,1}(X_{\pm})
\]
on the line bundle $\canclass_{X_{\pm}}^{n}$, given locally by
\[
\delbar_n=\frac{\del}{\del\zbar}\quad\text{and}\quad
\del_n=\rho^n\,\frac{\del}{\del z}\,\rho^{-n}.
\]
The $\delbar$-Laplacian acting on $\autoform_n(X_{\pm})$ is then
$\Delta_{n} = \delbar_n^*\delbar_n$,
where $\delbar_n^*=-\rho^{-1}\del_n$
is the adjoint of $\delbar_{n}$ with respect to
\eqref{metric-n-diff}.

The operator $\Delta_{n}$ is self-adjoint and non-negative, and has
pure discrete spectrum in the $L^{2}$-closure of
$\autoform_{n,m}(X_{\pm})$. The corresponding eigenvalues
$0\leq\lambda_0\leq\lambda_1\leq\dotsb$ of $\Delta_{n}$ have finite
multiplicity and accumulate only at infinity. The determinant of
$\Delta_n$ is defined by zeta regularization: the elliptic operator
zeta-function
\[
\zeta_{n}(s)=\sum_{\lambda_k>0}\lambda_k^{-s},
\]
defined initially for $\re s>1$, has a meromorphic continuation to the
entire $s$-plane, and by definition
\[
\det\Delta_n=e^{-\zeta_{n}'(0)}.
\]
The non-zero spectrum of $\Delta_{1-n}$ is identical to that of
$\Delta_{n}$, so that $\det\Delta_{n}=\det\Delta_{1-n}$. Hence
without loss of generality we will usually assume $n\geq 1$.

For $n\geq 2$, $\ker\Delta_{n}=\holoautoform_{n}(X_{\pm})$ has
dimension $d=(2n-1)(g-1)$.  If $\phi_{1},\dotsc,\phi_{d}$ is
a basis for  $\holoautoform_{n}(X_{\pm})$, we refer to
$(N_{n})_{jk}=\langle\phi_{j},\phi_{k}\rangle$ as the
\define{period matrix} corresponding to this basis.


\subsection{Quasi-Fuchsian deformation space}\label{QFdef}
(See \cite{Ahlfors2,Bers5,TT1}, and references therein.)
The set of marked, normalized quasi-Fuchsian groups of genus $g>1$
(recall our conventions on quasi-Fuchsian groups
from Section \ref{QFgroup}) has a natural structure of
a complex manifold of dimension $6g-6$.
We refer to it as the \define{quasi-Fuchsian space of genus $g$}
and denote it by $\qf_g$.  The subset of $\qf_{g}$ corresponding
to Fuchsian groups --- that is, the subset with
$X_{-} \simeq\ov{X_{+}}$ --- is called the Teichm\"uller space, $\teich_{g}$.
It is a totally real submanifold of $\qf_{g}$; however, it has a natural
complex structure.  With this complex structure, there is a natural
biholomorphism \cite{Kra2}
\begin{equation}\label{iso}
\qf_{g}=\teich_{g}\times\teich_{g}.
\end{equation}
We obtain the first (resp. second) copy of
$\teich_{g}$ by fixing $X_{-}$ (resp. $X_{+}$) (this is the
\define{Bers embedding} of $\teich_{g}$ into $\qf_{g}$).

Local coordinates (\define{Bers coordinates}) for $\qf_{g}$ may
be defined as follows.  Fix a
quasi-Fuchsian group $\Ga$ of genus $g$, which will serve as
the basepoint for the coordinate chart.  Let
$\harmform_{-1,1}(\Omega,\Ga)$ be the space of
\define{harmonic Beltrami differentials for $\Ga$}, that is,
$\mu\in\autoform_{-1,1}(\Omega,\Ga)$ such that $(\rho\mu)_{z}=0$. It
is a complex vector space of dimension $6g-6$; we give it the sup
norm.  Each $\mu$ in the open unit ball in
$\harmform_{-1,1}(\Omega,\Ga)$ is identified with a group $\Ga^\mu$,
a deformation of $\Ga$, by the generalized Riemann mapping theorem:
extending $\mu$ to $\widehat{\C}$ by $\mu\rvert_{\Lambda}=0$, there
exists a unique q.~c. homeomorphism $w_{\mu}: \widehat{\C}
\to\widehat{\C}$ satisfying the Beltrami equation
\begin{equation}\label{beltrami}
(w_{\mu})_{\z} = \mu (w_{\mu})_z
\end{equation}
and fixing the points $0$, $1$ and $\infty$. Set
\[
\Ga^{\mu} =w_{\mu}\Ga w_{\mu}^{-1}.
\]
The group $\Ga^{\mu}$ is quasi-Fuchsian, marked, and normalized, so it
may be identified with a point in $\qf_{g}$;
for a sufficiently small neighbourhood of the origin in
$\harmform_{-1,1}(\Omega,\Ga)$, this identification
is injective.

Write $U(\Ga)$ for the Bers coordinate chart based at $\Ga$.
There is a natural biholomorphism
$U(\Ga) \simeq U(\Ga^{\mu})$,
mapping $\Ga^{\nu} \in U(\Ga)$ to
$(\Ga^{\mu})^{\lambda} \in U(\Ga^{\mu})$ such that
$w_\nu=w_\lambda \circ w_\mu$,
which provides overlap maps for the coordinate charts, and allows
us to identify the holomorphic tangent space at the point
$\Ga^{\nu}\in U(\Ga)$ with the space of harmonic
Beltrami differentials $\harmform_{-1,1}(\Omega^{\nu},\Ga^{\nu})$.
Given $\mu \in \harmform_{-1,1}(\Omega,\Ga)$,
we denote by $\pa_{\mu}$ and $\pa_{\ov{\mu}}$ the holomorphic
and anti-holomorphic derivatives (vector fields) in a neighbourhood of
$\Ga$
defined using the Bers coordinates at the point $\Ga$.
The scalar product \eqref{metric-n-diff}
on $\harmform_{-1,1}(\Ga^{\nu})$ defines a K\"ahler metric on $\qf_{g}$
--- the \define{Weil-Petersson metric}.

To cover $\qf_{g}$, it is sufficient to take Bers coordinate charts
based at Fuchsian groups.  Explicitly: given
a Fuchsian group $\Ga_{0}$, let
$\mu_1,\ldots,\mu_{6g-6}$ be a real basis of
$\harmform_{-1,1}(\Ga_{0})$
which satisfies
\[
\mu(\z)=\ov{\mu(z)},
\]
and map $\boldsymbol{\vep}=(\vep_1, \ldots, \vep_{6g-6}) \mapsto
\Ga^{\boldsymbol{\vep}}= w_{\boldsymbol{\vep}}\circ \Ga_0\circ
w_{\boldsymbol{\vep}}^{-1}$, where
$w_{\boldsymbol{\vep}}:\hat{\C}\rightarrow \hat{\C}$ is the unique
normalized q.~c. mapping with Beltrami differential
$\vep_1\mu_1+\ldots\vep_{6g-6}\mu_{6g-6}$.  In this coordinate
chart, the subset $\im\vep_{k}=0$, $k=1,\dotsc,6g-6$ consists of
precisely the Fuchsian groups in the chart.  These coordinate charts
cover $\qf_{g}$.

For $\Ga$ quasi-Fuchsian and $\mu$ in the unit ball in
$\harmform_{-1,1}(\Ga)$, if $w_{\vep}$ is the corresponding q.~c.
map satisfying $\eqref{beltrami}$, we then have, for each
$z\in\Omega$,
\begin{equation}\label{var2}
\left.\frac{\pa}{\pa \ov{\vep}}\right|_{\vep=0} w_{\vep\mu}(z) = 0.
\end{equation}
Combining this with the equation
\begin{equation}\label{group}
\gamma^{\vep\mu} \circ w_{\vep\mu}=w_{\vep\mu}\circ \gamma,
\end{equation}
we have
\begin{equation}\label{var3}
\frac{\pa}{\pa\ov{\vep}}\Bigr|_{\vep=0} \gamma^{\vep\mu}=0
\end{equation}
for each $\gamma^{\vep\mu}\in\Ga^{\vep\mu}$.


\subsection{Families of \n differentials}
(See \cite{Bers5} and references therein.)
The quasi-Fuchsian fibre space is a fibration
$p:\univqf_{g}\to\qf_{g}$ with fibre
$\pi^{-1}(t)=\Ga^{t}\bk\Omega^{t}$
for $t\in\qf_{g}$.  Let
$T_V\univqf_g\to\univqf_g$ be  the holomorphic vertical
tangent bundle --- the  holomorphic line bundle over $\univqf_g$
consisting of vectors in the holomorphic tangent space
$T\univqf_g$ that are tangent to the fibres $\pi^{-1}(t)$.
A  family $\phi^{t}$ of \n differentials
is defined as a smooth section of the direct image bundle
\[
\Lambda_{n}=p_{*}((T_V\univqf_g)^{-n})\to\qf_{g}.
\]
The fibre of $\Lambda_{n}$ over $t\in\qf_{g}$ is the vector space $\autoform_n(\Omega^{t},\Ga^{t})$.
The hyperbolic metric $\rho^{t}$ on $\Ga^{t}\bk\Omega^{t}$
defines a natural Hermitian metric on the line bundle
$\Lambda_{n}$ by \eqref{metric-n-diff}.

Analogously, there exist fibre spaces $p_{\pm}:(\univqf_{g})_{\pm}\to\qf_{g}$ with fibre
$\pi^{-1}(t)=\Ga^{t}\bk\Omega_{\pm}^{t}$
for $t\in\qf_{g}$, and consequently we may define families of
\n differentials in $\autoform_n(\Omega_{\pm}^t,\Gamma^t)$.

For a harmonic Beltrami differential $\mu$ in the unit ball of
$\harmform_{-1,1}(\Omega,\Ga)$, the pullback of an \n differential
$\phi^{\vep}\in\autoform_n(\Omega^{\vep\mu},\Ga^{\vep\mu})$ is an \n
differential $w^*_{\vep\mu}(\phi^\vep)\in\autoform_n(\Omega,\Ga)$
defined by
\begin{equation*}
w^*_{\vep\mu}(\phi^\vep)
=\phi^\vep\circ w_{\vep\mu}
\,\bigl((w_{\vep\mu})_{z}\bigr)^n,
\end{equation*}
where $w_{\vep\mu}:\widehat\C\to\widehat\C$ is the
solution of the Beltrami equation corresponding to $\mu$.
The Lie derivatives of the
family $\phi^\vep$ in the directions $\mu$ and
$\overline{\mu}$ are defined by
\begin{align*}
L_\mu\phi
&=\frac{\del}{\del\vep}\bigg\vert_{\vep=0}
 w^*_{\vep\mu}(\phi^\vep) \in\autoform_n(\Omega,\Ga)
\\
\text{and}\quad
L_{\ov\mu}\phi
&=\frac{\del}{\del\overline{\vep}}\bigg\vert_{\vep=0}
 w^*_{\vep\mu}(\phi^\vep) \in\autoform_n(\Omega,\Ga).
\end{align*}
Note that a smooth section $\phi^{t}$ of $(\Lambda_{n})_{\pm}$ is
holomorphic if and only if $L_{\ov{\mu}}\phi=0$ in each Bers
coordinate chart.  If $\phi^t_1,\dotsc,\phi^t_d$ are holomorphic
sections which form a basis of
$\holoautoform_{n}(\Omega_{\pm}^{t},\Ga^{t})$ for each
$t\in\qf_{g}$, we call them a \define{holomorphically varying basis}
of holomorphic \n differentials.


\subsection{Inversion on the quasi-Fuchsian deformation space}
\label{inversion}
There is a canonical inversion $\iota$ on
the quasi-Fuchsian deformation space given by
$\Ga \mapsto\ov{\Ga}$, where $\ov{\Ga}$
is the quasi-Fuchsian group
\begin{align}\label{bargroup}
\ov{\Ga}=\left\{
\ov{\gamma}=\begin{pmatrix}\ov{a}&\ov{b}\\\ov{c}&\ov{d}\end{pmatrix}\;;\;
\gamma=\begin{pmatrix} a&b\\c&d\end{pmatrix}\in \Ga\right\}.
\end{align}
With our conventions, we have
$\Omega_{\pm}(\ov{\Ga})=\ov{\Omega_{\mp}(\Gamma)}$,
where $\ov{\Omega_{\pm}}$ means the set of $z$ such that
$\ov{z}\in\Omega_{\pm}$, with the orientation from the complex plane.
The hyperbolic metric on $\ov{\Omega}$ is given by
$\rho_{\ov{\Omega}}(z)=\rho_{\Omega}(\ov{z})$.
If $\Ga\bk\Omega\simeq X_+\sqcup X_-$, then
$\ov{\Ga}\bk \ov{\Omega} \simeq \ov{X_-}\sqcup \ov{X_+}$,
where $\ov{X_{\pm}}$ is obtained from $X_{\pm}$ by taking the conjugate
of all local coordinates.  Note that as a hyperbolic manifold, $\ov{X_{\pm}}$ is
isometric to $X_{\pm}$ by an orientation-reversing isometry.  The
fixed submanifold of the inversion $\iota$ is precisely the
real submanifold of Fuchsian groups.   On the Bers coordinate with
a Fuchsian basepoint $\Ga_{0}$, the inversion is realized
explicitly by $\iota(\mu)(z)=\ov{\mu(\ov{z})}$.


\subsection{Classical Liouville action}\label{action} In \cite{TT1}, Takhtajan and
the second author constructed the classical Liouville action
$S:\qf_g\to \R$ on the quasi-Fuchsian
space, which has the following properties:
\begin{itemize}
\item[\textbf{CL1}] $S$ is a real analytic function on
$\qf_g$.

\item[\textbf{CL2}] $S$ is invariant under the inversion
$\iota$, that is, $S(\Ga)= S(\ov{\Ga})$.

\item[\textbf{CL3}] Restricted to the real submanifold of Fuchsian
groups, $S$ is a constant, equal to $8\pi(2g-2)$.

\item[\textbf{CL4}] $-S$ is a potential for the
Weil-Petersson metric on the quasi-Fuchsian deformation space.
Namely, at a point $\Ga$, for all  $\mu, \nu \in \Omega_{-1,1}(\Ga)$,
\[
\pa_{\mu}\pa_{\ov{\nu}}(-S) = \langle \mu, \nu\rangle.
\]

\end{itemize}


\subsection{Selberg Zeta function and D'Hoker-Phong theorem}
The Selberg zeta function $Z(s)$ of a Riemann surface $X
\simeq\Ga\bk \hp_+$, where $\Ga$ is a Fuchsian group, is defined for
$\text{Re}\; s>1$ by the absolutely convergent product
\[
Z(s)
= \prod_{\{\g_0\}} \prod_{m=0}^{\infty} ( 1 -\lambda(\g_0)^{s+m}),
\]
where $\g_0 $ runs over the set of conjugacy classes of primitive
hyperbolic elements of $\Ga$, and $0<\lambda(\gamma) < 1$ is the
multiplier of $\gamma$. The function $Z(s)$ has a meromorphic
continuation to the whole $s$-plane. In \cite{DhokerPhong},
D'Hoker and Phong proved that for $n\geq 2$,
\begin{align}\label{DP}
\det \Delta_n = c_{g,n} Z(n).
\end{align}
Here $c_{g,n}$ is  an explicitly known
positive constant depending only on $g$ and $n$, and not on $X$
(we refer the reader to \cite{DhokerPhong} for the precise
expression for $c_{g,n}$).

\section{Generalized local index theorem on the quasi-Fuchsian space}
\label{localindex}

In \cite{ZT} (see also \cite{TZ91}), Takhtajan
and Zograf proved the Local Index Theorem:
\begin{theorem*}
Let $\Ga$ be a marked Fuchsian group corresponding to a point in the
Teichm\"uller space $\teich_{g}$, with $\Ga\bk\hp_+\simeq X$ a
compact Riemann surface of genus $g>1$.  Write $[X]$ for the
corresponding marked Riemann surface.  Then we have, for all
$\mu,\nu\in \Omega_{-1,1}(\hp_+,\Ga)$,
\[
\pa_{\mu}\pa_{\ov{\nu}}\log\frac{\det\Delta_n(X)}{\det N_{n}([X])}
= \frac{6n^2-6n+1}{12\pi}\langle \mu,\nu\rangle,
\]
where $(N_n([X]))_{jk}=\langle\phi_{j},\phi_{k}\rangle$ is the
period matrix of a holomorphically varying basis
$\phi_{1},\dotsc,\phi_{d}$ of $\holoautoform_{n}(X)$, i.~e. of holomorphic
n-differentials of the Riemann surface $X$.
\end{theorem*}
Note that we write $N_{n}([X])$ since we may find a holomorphically varying
basis globally only over $\teich_{g}$, not over the moduli space of surfaces $X$.

Given $\Ga$ a quasi-Fuchsian group with
$\Ga\bk\Omega\simeq X\simeq X_{+}\sqcup X_{-}$
a union of compact Riemann surfaces $X_\pm$ of genus $g>1$, the
determinant of the Laplacian acting on \n differentials
is the product
\[
\det\Delta_n(X)
=\det\Delta_n(X_{+})\det\Delta_n(X_{-})
=c_{g,n}^2Z_{X_+}(n)Z_{X_-}(n),
\]
and, given choices of holomorphically varying bases of holomorphic
\n differentials on $X_+$ and $X_-$, the period matrix $N_n([X])$ is
defined as
\[
\det N_n([X])
=\det N_{n}([X_{+}])\det N_{n}([X_{-}]).
\]

It is straightforward to generalize the Local Index Theorem to
quasi-Fuchsian groups:
\begin{theorem*}
Let $\Ga$ be a marked quasi-Fuchsian group corresponding to a point
in the quasi-Fuchsian space $\qf_{g}$, with $\Ga\bk\Omega\simeq
X\simeq X_{+}\sqcup X_{-}$ a union of two compact Riemann surfaces
of genus $g>1$. Write $[X_{\pm}]$ for the corresponding marked
Riemann surfaces.  Then we have, for all $\mu,\nu\in
\Omega_{-1,1}(\Omega,\Ga)$,
\[
\pa_{\mu}\pa_{\ov{\nu}}\log\frac{\det\Delta_n(X)}{\det N_{n}([X])}
=\pa_{\mu}\pa_{\ov{\nu}}
 \log\frac{\det\Delta_n(X_+)}{\det N_n([X_+])}
   \frac{\det\Delta_n(X_-)}{\det N_n([X_-])}
=\frac{6n^2-6n+1}{12\pi}\langle \mu,\nu\rangle.
\]
We call this the Generalized Local Index Theorem.
\end{theorem*}
\begin{proof}
Simply use the isomorphism
$\qf_g\simeq\teich_g\times\teich_g$ to split the tangent
space $T_t\qf_g\simeq T_t\teich_g\oplus T_t\teich_g$
at each point $t\in\qf_g$.
\end{proof}

\section{The function $F(n)$}\label{Ffunction}

Analogous to the Selberg zeta function $Z(s)$, given a quasi-Fuchsian
group $\Ga$ and an integer $n\geq 2$, we define the function
\[
F(n) = \prod_{\{\g\}} \prod_{m=0}^{\infty}
  ( 1 -\lambda(\g_0)^{n+m}),
\]
where $\{\g\} $ runs over the set of conjugacy classes of primitive
elements of $\Ga$, omitting the identity, and $\tabs{\lambda(\g)}< 1$ the
multiplier of $\g$.  The product
converges absolutely if and only if the series
$\sum_{\{\gamma\}}\sum_{m=0}^{\infty}\tabs{\lambda_{\gamma}}^{m+n}$
converges.  It is straightforward to prove that this series converges provided that the
multiplier series $\sum_{[\gamma]}\tabs{\lambda_{\gamma}}^{n}$ converges,
where $[\gamma]$ runs over all distinct conjugacy classes (not necessarily primitive) in $\Ga$.  The argument of B\"user
for Schottky groups \cite{Buser} goes through for
quasi-Fuchsian groups, showing that multiplier series converges if the
Poincar\'e series $\sum_{[\gamma]}\tabs{\g'(z)}^{n}$ converges, and it is a
classical fact that this converges when $n\geq 2$.
\begin{lemma}\label{holF}\mbox{}
\begin{itemize}
\item[(i)] $F(n)$ is a holomorphic function on the
quasi-Fuchsian deformation space.
\item[(ii)] Under the inversion $\iota$ on the quasi-Fuchsian deformation
space, the function $F(n)$ transforms as
$F(n)(\ov{\Ga})=\ov{F(n)(\Ga)}$.
\item[(iii)] Restricted to the real submanifold of Fuchsian groups, $F(n)$ is real and
coincides with $Z(n)$.
\end{itemize}
\end{lemma}
\begin{proof}
Properties (ii) and (iii) are immediate from the definitions.  It
is easy to show that the multiplier $\lambda(\gamma)$
is a holomorphic function of the entries of $\gamma$, whenever
$0<\lambda(\gamma)<1$.  Combining this with
\eqref{var3} establishes property (i).
\end{proof}
These properties characterize $F(n)$ uniquely; see Lemma
\ref{reflection}.

\section{Bers integral operator and Bers dual bases}

Let $\Ga$ be a quasi-Fuchsian group of genus $g>1$
with ordinary set $\Omega_{+}\sqcup\Omega_{-}$
(recall our conventions on quasi-Fuchsian groups
from Section \ref{QFgroup}),
and let $n\geq 2$ be an integer. Recall that
$A_n(\Omega_\pm, \Ga)$, the space of holomorphic
\n differentials of $\Ga$ with support on $\Omega_\pm$,
is a complex vector space of dimension
$d= (2n-1)(g-1)$,  with a canonical inner product given by
\[\label{inner}
\langle \phi, \psi\rangle
=\int\limits_{\Ga\bk\Omega_\pm} \phi(z)\ov{\psi(z)}\rho(z)^{1-n}~\ud^2 z.
\]


\subsection{Bers integral operator}
In \cite{Bers2}, Bers introduced complex linear operators
\[
K_\pm(\Ga):\ov{A_n(\Omega_\pm,\Ga)}\to A_n(\Omega_\mp,\Ga)
\]
(we suppress the $n$ to simplify notation).
The two operators $K_\pm(\Ga)$ (or simply $K_\pm$) are defined for
$\phi\in A_n(\Omega_\pm,\Ga)$ and $z\in\Omega_\mp$ by
\[
(K_\pm\ov{\phi})(z)
=c_n\int\limits_{\Omega_\pm}\frac{\ov{\phi(w)}\rho(w)^{1-n}}{(z-w)^{2n}}~\ud^2w
=\int\limits_{\Ga\bk\Omega_\pm}K_\pm(z,w)\ov{\phi(w)}\rho(w)^{1-n}~\ud^2w,
\]
where $c_n$ is the constant
\[
c_n=\frac{2^{2n-2}(2n-1)}{\pi},
\]
and, for $z\in\Omega_\mp$, $w\in\Omega_\pm$,
\[
K_\pm(z,w)
=c_n\sum_{\gamma\in\Ga}
\frac{\gamma'(z)^n}{(\gamma z-w)^{2n}}.
\]
We can also define these operators for the conjugate group
$\ov{\Ga}$; we have
\[
K_\pm(\ov{\Ga}):\ov{A_n(\ov{\Omega_\mp},\ov{\Ga})}\to A_n(\ov{\Omega_\pm},\ov{\Ga}).
\]

\begin{remark}
In \cite{Bers2}, the operator $K_-$ is defined on the Banach space
$\ov{\mathcal{B}_n(\Omega_-, \Ga)}$, where
$\mathcal{B}_n(\Omega_-, \Ga)$ is the space of bounded
\n differentials of $\Ga$ with support on $\Omega_-$. There is a
canonical decomposition
\[
\mathcal{B}_n = A_n \oplus \mathcal{N}_n,
\]
where $\mathcal{N}_n$
is the subspace of $\mathcal{B}_n$ consisting of all $\eta$ such
that
\[
\int\limits_{\Ga\bk\Omega_-}\eta\ov{\phi}\rho^{1-n}=0,
\]
for all holomorphic \n differentials $\phi$ which are $L^1$ integrable.
Using this characterization of $\mathcal{N}_n$, it is easy to see
that $\ov{\mathcal{N}_n}$ lies in the kernel of $K_-$. Hence here
we define $K_-$ on the quotient space
$\ov{\mathcal{B}_n}/\ov{\mathcal{N}_n}\simeq \ov{A_n}$
(and similarly for $K_+$).
\end{remark}

Define operators
\begin{align*}
\iota_\pm&:A_n(\Omega_\pm,\Ga)\to A_n(\ov{\Omega_\pm},\ov{\Ga}),\\
\ov{\iota_\pm}&:\ov{A_n(\Omega_\pm,\Ga)}\to \ov{A_n(\ov{\Omega_\pm},\ov{\Ga})}
\end{align*}
by $(\iota_\pm\phi)(z)=\ov{\phi(\ov{z})}$ and $(\ov{\iota_\pm}\,\ov{\phi})(z)=\phi(\ov{z})$
for each $\phi\in A_n(\Omega_\pm,\Ga)$ and $z\in\ov{\Omega_\pm}$.

\begin{lemma}\label{propKL1} \mbox{}
\begin{itemize}
\item[(i)] $K_\pm^*\phi=\ov{K_\mp\ov{\phi}}$
for all $\phi$ in $A_{n}(\Omega_\mp,\Ga)$.
\item[(ii)] $\iota_\mp K_\pm(\Ga)=K_\mp(\ov{\Ga})\ov{\iota_\pm}$ and
$\ov{\iota_\pm}K^*_\pm(\Ga)=K^*_\mp(\ov{\Ga})\iota_\mp$.
\item[(iii)] When $\Ga$ is Fuchsian,
$K_\pm\ov{\phi}=\iota_\pm\phi$ for all $\phi$ in $A_n(\Omega_\pm,\Ga)$.
\end{itemize}
\end{lemma}
\begin{proof}
For $z\in\Omega_\pm$ and $\phi\in A_{n}(\Omega_\mp,\Ga)$, the adjoint operator
$K_\pm^*:A_n(\Omega_\mp,\Ga)\rightarrow\ov{A_n(\Omega_\pm,\Ga)}$
is given by
\[
(K_\pm^*\phi)(z)
=c_n\int\limits_{\Omega_\mp}\frac{\phi(w)\rho(w)^{1-n}}{(\zbar-\wbar)^{2n}}~\ud^2w,
\]
hence
\[
\ov{K_\pm^*\phi(z)}
=c_n\int_{\Omega_\mp}\frac{\ov{\phi(w)}\rho(w)^{1-n}}{(z-w)^{2n}}~\ud^{2}w
=(K_{\mp}\ov{\phi})(z),
\]
establishing (i).  For $z\in\ov{\Omega_\mp}$ and $\phi\in A_{n}(\Omega_\pm,\Ga)$,
\[
\begin{split}
(\iota_\mp K_\pm(\Gamma)\ov{\phi})(z)
=\ov{(K_\pm(\Ga)\ov{\phi})(\zbar)}
&=c_{n}\int_{\Omega_\pm}\frac{\phi(w)\rho(w)^{1-n}}{(z-\wbar)^{2n}}~\ud^{2}w \\
&=c_{n}\int_{\ov{\Omega_\pm}}\frac{\phi(\wbar)\rho(\wbar)^{1-n}}{(z-w)^{2n}}~\ud^{2}w\\
&=(K_\mp(\ov{\Ga})\ov{\iota_\pm}\,\ov{\phi})(z),
\end{split}
\]
which proves the first part of (ii); the second part is similar.
When $\Ga$ is Fuchsian, property (iii) is the reproducing formula
--- see \cite{Bers2}.
\end{proof}

Bers proved that the operator $K_\pm$ is invertible.
Hence $\kappa_\pm=K_\pm K_\pm^*$ is a self adjoint positive
definite operator.
\begin{lemma}\label{propKL2}\mbox{}
\begin{itemize}
\item[(i)] $\kappa_+$ and $\kappa_-$ have the same eigenvalues.
\item[(ii)] $\kappa_\pm(\ov{\Ga})$ and $\kappa_\mp(\Ga)$ have the same eigenvalues.
\item[(iii)] If $\Ga$ is Fuchsian, $\kappa_\pm(\Ga)$ is the identity.
\end{itemize}
\end{lemma}
Note that properties (i) and (ii) imply that $\kappa_\pm(\Ga)$ and
$\kappa_\pm(\ov{\Ga})$ have the same eigenvalues.
\begin{proof}
The map $\phi\mapsto K_+\ov{\phi}$ from
$A_n(\Omega_+,\Ga)\to A_n(\Omega_-,\Ga)$ conjugates
$\kappa_+$ with $\kappa_-$:
\[
K_+K_+^*(K_+\ov{\phi})
=K_+K_+^*\ov{K_-^*\phi}
=K_+\ov{K_-K_-^*\phi},
\]
which establishes (i).  To prove (ii), note that the map $\iota_-$ conjugates
$\kappa_+(\ov{\Ga})$ with $\kappa_-(\Ga)$:
\[
\iota_-K_+(\Ga)K_+^*(\Ga)
=K_-(\ov{\Ga})\ov{\iota_+}K_+^*(\Ga)
=K_-(\ov{\Ga})K_-^*(\ov{\Ga})\iota_-.
\]
Property (iii) follows from part (iii) of the previous lemma.
\end{proof}


\subsection{Bers dual bases}\label{period}
We choose a basis $\phi_1^+,\ldots,\phi_d^+$ of $A_n(\Omega_+,\Ga)$
and expand the kernel $K_-(z,w)$ with respect to this basis by
\[\label{basis}
K_-(z,w)=\sum_{k=1}^d \phi_k^+(z)\phi_k^-(w) .
\]
We define the period matrices $N_\pm=N_n([X_\pm])$ by
\[
(N_\pm)_{kl} = \langle \phi_k^\pm,\phi_l^\pm\rangle.
\]
Then
\begin{align*}
\kappa_-(z,w)
=\int\limits_{\Ga\bk\Omega_-}K_-(z,u)\ov{K_-(w,u)}\rho(u)^{1-n}~\ud^2u
= \sum_{j =1}^d\sum_{k=1}^d\phi_j^+(z) (N_-)_{jk} \ov{\phi_k^+(w)},
\end{align*}
and hence
\[
(\kappa_-\phi_l^+)(z)
= \sum_{j=1}^d \left(\sum_{k=1}^d
  (N_+)_{lk}(N_-)_{jk}\right) \phi_j^+(z).
\]
Namely, with respect to the basis $\phi_1^+,\ldots, \phi_d^+$,
the matrix for the operator $\kappa_-$ is given by $N_+N_-^T$. The
invertibility of $\kappa_-$ then shows that
$\phi_1^-, \ldots,\phi_d^-$ is a basis of $A_n(\Omega_-,\Ga)$,
and we say that
this basis is \define{Bers dual} to $\phi_1^+, \ldots, \phi_d^+$.
Note that we have
$\det\kappa_-=\det N_+\det N_-$.


\subsection{Properties of period matrices}\label{holbasis}
Suppose that
$\phi_1^+,\dotsc,\phi_d^+$ is a basis
for $A_{n}(\Omega_+,\Ga)$ chosen globally on
$\qf_g$ and varying holomorphically, and
$\phi_1^-,\dotsc,\phi_d^-$ is a Bers dual basis
for $A_{n}(\Omega_-,\Ga)$.  Let $N_\pm$
be the corresponding period matrices.
\begin{lemma}\label{invK}\mbox{}
\begin{itemize}
\item[(i)] The basis $\phi_1^-,\dotsc,\phi_d^-$ varies holomorphically on $\qf_g$.
\item[(ii)] $\det N_+\det N_-$ is invariant under the inversion $\iota$ on $\qf_g$.
\item[(iii)] $\det N_+\det N_- =1$ on the real submanifold of Fuchsian groups.
\end{itemize}
\end{lemma}
\begin{proof}
If $\mu\in\Omega_{-1,1}(\Ga)$ is a harmonic Beltrami differential, and
$w_{\vep\mu}$ is the q.~c. mapping defined by
$\vep\mu$, then by \eqref{var2} and \eqref{group}, we have
\[
\frac{\pa}{\pa\ov{\vep}}\Bigr|_{\vep=0} \left((w_{\vep\mu}\boxtimes
w_{\vep\mu} )^*K_-\right)(z,w)= 0,
\]
the symbol $\boxtimes$ indicating that we are pulling $K_-(z,w)$ back as
an \n differential in each variable.
Hence,
\begin{align*}
\sum_{k=1}^{d}
 \left(\frac{\pa}{\pa\ov{\vep}}\Bigr|_{\vep=0}
   w_{\vep\mu}^*(\phi_k^-)^{\vep\mu}\right)(w)(\phi_k^+)^{\vep\mu}(z)
+ (\phi_k^-)^{\vep\mu}(w)
 \left(\frac{\pa}{\pa\ov{\vep}}\Bigr|_{\vep=0}
   w_{\vep\mu}^*(\phi_k^+)^{\vep\mu}\right)(z)=0.
\end{align*}
Consequently, if the basis $\phi_1^+,\ldots, \phi_d^+$
varies holomorphically with respect to moduli, i.~e. for all $k$,
\[
\left(\frac{\pa}{\pa\ov{\vep}}\Bigr|_{\vep=0}
 w_{\vep\mu}^*(\phi_k^+)^{\vep\mu}\right)(z)=0,
\]
then
\[
\left(\frac{\pa}{\pa\ov{\vep}}\Bigr|_{\vep=0}
 w_{\vep\mu}^*(\phi_k^-)^{\vep\mu}\right)(z)=0
\]
for all $k$ as well, which proves (i).

Properties (ii) and (iii) follows from Lemma \ref{propKL2}, since
$\det\kappa_-=\det N_+\det N_-$.
\end{proof}

\subsection{Global choice of basis over $\qf_g$}
We may choose a holomorphically varying basis of \n differentials
globally over the Teichm\"uller space by using the Bers embedding.
Fix a marked Riemann surface $[X_+]$, and choose a fixed basis
$\phi_1^+,\dotsc,\phi_d^+$ of $\holoautoform_n(X_+)$. Then for any
$[X_-]$, take the corresponding quasi-Fuchsian group, identify
$\phi_1^+,\dotsc,\phi_d^+$ with a basis of
$\holoautoform_n(\Omega_+,\Ga)$, and let $\phi_1^-,\dotsc,\phi_d^-$
be the Bers dual basis of $\holoautoform_n(\Omega_-,\Ga)$.  By Lemma
\ref{invK}, this basis varies holomorphically. (This is essentially
equivalent to a construction by Bers; see \cite{Bers5} and
references therein.)

Now, assume a choice of holomorphically varying basis of \n
differentials $\phi_1^+,\dotsc,\phi_d^+$ of $\holoautoform_n(X_+)$
has been made at each point $[X_+]$ in the Teichm\"uller space.
Given any $[X_-]$, take the corresponding quasi-Fuchsian group,
identify $\phi_1^+,\dotsc,\phi_d^+$ with a basis of
$\holoautoform_n(\Omega_+,\Ga)$, and let $\phi_1^-,\dotsc,\phi_d^-$
be the Bers dual basis of $\holoautoform_n(\Omega_-,\Ga)$.  By Lemma
\ref{invK}, this basis varies holomorphically.

Consequently we have shown that it is possible to make a choice
of holomorphically varying bases
$\phi_1^\pm,\dotsc,\phi_d^\pm$
of $\holoautoform_n(\Omega_\pm,\Ga)$ globally over
$\qf_g$, in such a way that $\phi_1^-,\dotsc,\phi_d^-$
is Bers dual to $\phi_1^+,\dotsc,\phi_d^+$ at each point of $\qf_g$.

\section{Holomorphic factorization of determinants of Laplacians}\label{mainthm}

In this section, we prove our main theorem:
\begin{theorem*}
Let $\Ga$ be a quasi-Fuchsian group simultaneously uniformizing
compact Riemann surfaces $X_+$ and $X_-$ of genus $g>1$. Then for
$n\geq 2$,
\[
\frac{\det\Delta_n(X_+)}{\det N_n([X_+])}\frac{\det\Delta_n(X_-)}{\det N_n([X-])}
= a_{g,n}\abs{F_{\Ga}(n)}^2
\exp\left(-\frac{6n^2-6n+1}{12\pi}S_{\Ga}\right).
\]
Here the Laplacian is computed in the hyperbolic metric;
$N_n([X_\pm])$ are period matrices of bases
of  $\ker\Delta_{n}(X_\pm)$,
Bers dual in the sense of Section \ref{period}
and chosen globally over the quasi-Fuchsian space;
and $S_{\Ga}$ is the classical Liouville action defined in
Section \ref{action}.
The function
\[
F_{\Ga}(n)
=\prod_{\{\gamma\}\in\Ga}\prod_{m=0}^{\infty}
\big(1-\lambda(\gamma)^{n+m}\big)
\]
is defined in more detail in Section \ref{Ffunction};
and $a_{g,n}$ is the positive constant
\[
a_{g,n}=c_{g,n}^2 \exp\left( \frac{(6n^2-6n+1)(4g-4)}{3}\right),
\]
where $c_{g,n}$ is the constant from the D'Hoker-Phong formula
\eqref{DP}.
\end{theorem*}
First, we prove a ``Schwarz reflection'' lemma for pluriharmonic functions:
\begin{lemma}\label{reflection}
Let $V$ be an open, convex subset of $\CC^m$, such that
$V\cap\RR^m$ is nonempty.  Suppose that $h:V\to\RR$ satisfies
the following conditions:
\begin{itemize}
\item[(i)] $\displaystyle{\frac{\partial^{2}h}{\partial z_{j}\partial\overline{z}_{k}}=0}$
everywhere in $V$, for all $j,k\in\{1,\dotsc,m\}$.
\item[(ii)] $h(x_1,\dotsc,x_m)=0$ for all $(x_1,\dotsc,x_m)\in V\cap\RR^{m}$.
\item[(iii)] $h(\overline{z}_1,\dotsc,\overline{z}_m)=h(z_1,\dotsc,z_m)$
for all $(z_1,\dotsc,z_m)\in V$
\\ \hspace*{0cm} such that
$(\overline{z}_1,\dotsc,\overline{z}_m)\in V$.
\end{itemize}
Then $h$ is identically zero on $V$.
\end{lemma}
\begin{proof}
Fix $a,b\in\RR^m$, and let $W=\{z\in\CC:az+b\in V\}$.  Define
$f:W\to\RR$ by $f(z)=h(az+b)$.  Property (i) implies that
$f$ is harmonic on $W$. Property (ii) implies that $f(x)=0$
for all $x\in W\cap\RR$. Hence by the Schwarz reflection
principle, we have $f(\ov{z})=-f(z)$ for all $z$ in some
neighbourhood of $W\cap\RR$.  On the other hand, $f(\ov{z})=f(z)$
whenever $z,\ov{z}\in W$ by (iii), so $f(z)=0$ for all $z$
in some neighbourhood of $W\cap\RR$.  But since $V$ is convex, $W$
is connected, so $f(z)=0$ for all $z\in W$.

Hence $h$ is zero at all points in $V$ of the form $az+b$ for some
$a,b\in\RR^m$ and some $z\in\CC$.  But it is easy to check that
every point in $\CC^m$ is of this form, so $h$ is identically
zero on $V$.
\end{proof}

Now we return to the proof of the theorem.  Let
\[
h= \log\left(\left.
\frac{\det\Delta_n (X_+)}{\det N_n([X_+])}\frac{ \det\Delta_n(X_-)}{\det N_n([X_-])}
\right/
a_{g,n} \abs{F_{\Ga}(n)}^2
\exp\left(-\frac{6n^2-6n+1}{12\pi}S_{\Ga}\right)
\right).
\]
The function $h$ is real-valued on the quasi-Fuchsian deformation
space $\qf_g$.  We claim that
 \begin{itemize}
 \item[(i)] $h$ is pluriharmonic on $\qf_g$
 \item[(ii)] $h$ is invariant under the inversion $\iota$
 \item[(iii)] $h=0$ on the real submanifold of Fuchsian groups.
 \end{itemize}
Property (i) follows directly from the generalized local index theorem
(Section \ref{localindex}) and the fact that $F_{\Ga}(n)$ is holomorphic on
$\qf_g$ (Section \ref{Ffunction}).  Property (ii) follows from the transformations
under $\iota$, established in Sections \ref{detLaplacian},
\ref{Ffunction}, \ref{action}, and \ref{holbasis},
of the factors appearing in
$h$.  Using the fact that $S_{\Ga}$ is constant on the real submanifold of
Fuchsian groups, property (iii) reduces to the D'Hoker-Phong formula \eqref{DP}.

Writing $h$ in the local coordinates on $\qf_g$ described in Section
\ref{QFdef}, we see that it satisfies the conditions of Lemma
\ref{reflection}, and hence is identically zero on $\qf_g$, proving
the theorem.

\begin{remark}
The theorem should be considered as an equality of functions
over the quasi-Fuchsian space $\qf_g$.  However, we can obtain as
a corollary an equality of functions over the Teichm\"uller space $\teich_g$,
by means of the Bers embedding.  Fix the marked Riemann surface
$[X_-]$ and the basis $\phi_1^-,\dotsc,\phi_d^-$ of holomorphic
\n differentials on $[X_-]$.  Then
$\det\Delta_n(X_-)/\det N_n([X_-])$ is constant, and all other quantities
in the theorem are functions only of $[X_+]\in\teich_g$.
\end{remark}


\end{document}